\documentclass{gtart}

\input gtoutput
\volumenumber{3}\papernumber{15}\volumeyear{1999}
\pagenumbers{369}{396}\published{25 October 1999}
\proposed{Ralph Cohen}
\seconded{Robion Kirby, Steve Ferry}
\received{28 June 1999}
\accepted{21 October 1999}

\usepackage{amssymb, amsmath, amscd}

\newtheorem{theorem}{Theorem}[section]

\newtheorem{lemma}[theorem]{Lemma}

\theoremstyle{definition}
\newtheorem{notation}[theorem]{Notation}
\newtheorem{definition}[theorem]{Definition}
\newtheorem{remark}[theorem]{Remark}
\newtheorem{remarks}[theorem]{Remarks}
\newtheorem{example}[theorem]{Example}
\newtheorem{assumption}[theorem]{Assumption}
\newtheorem{conjecture}[theorem]{Conjecture}

\newtheorem*{update}{Update}
\newtheorem*{acknowledgments}{Acknowledgments}

\newcommand{\eqdef}{\;{:=}\;}
\newcommand{\C}{{\mathbb C}}
\newcommand{\Q}{{\mathbb Q}}
\newcommand{\R}{{\mathbb R}}
\newcommand{\Z}{{\mathbb Z}}
\newcommand{\op}{\operatorname}
\newcommand{\mc}[1]{{\mathcal #1}}
\newcommand{\Spinc}{\op{Spin}^c}

\newcommand{\Ker}{\op{Ker}}

\newcommand{\SW}{\op{SW}}
\newcommand{\Nov}{\op{Nov}}
\newcommand{\tensor}{\otimes}
\newcommand{\Gr}{\op{Gr}}

\newcommand{\Cnov}{CN}
\newcommand{\Eul}{\op{Eul}}
\def\S{section }

\begin{document}

\title{Circle-valued Morse theory and Reidemeister torsion}                    
\authors{Michael Hutchings\\Yi-Jen Lee}                  
\address{Dept of Math, Stanford University, Stanford, CA 94305, USA\\
\smallskip\\
Dept of Math, Princeton University, Princeton, NJ 08544, USA}
\asciiaddress{Dept of Math, Stanford University\\Stanford, CA 94305, USA\\
\\Dept of Math, Princeton University\\Princeton, NJ 08544, USA}
\email{hutching@math.stanford.edu\\ylee@math.princeton.edu}

\asciiabstract{Let X be a closed manifold with zero Euler characteristic, 
and let f: X --> S^1 be a circle-valued Morse function.  We define an
invariant I which counts closed orbits of the gradient of f, together
with flow lines between the critical points.  We show that our
invariant equals a form of topological Reidemeister torsion defined by
Turaev [Math. Res. Lett. 4 (1997) 679-695].

We proved a similar result in our previous paper [Topology, 38 (1999)
861-888], but the present paper refines this by separating closed
orbits and flow lines according to their homology classes.
(Previously we only considered their intersection numbers with a fixed
level set.)  The proof here is independent of the previous proof and
also simpler.

Aside from its Morse-theoretic interest, this work is motivated by the
fact that when X is three-dimensional and b_1(X)>0, the invariant I
equals a counting invariant I_3(X) which was conjectured in our
previous paper to equal the Seiberg-Witten invariant of X.  Our
result, together with this conjecture, implies that the Seiberg-Witten
invariant equals the Turaev torsion.  This was conjectured by Turaev
[Math. Res. Lett. 4 (1997) 679-695] and refines the theorem of Meng
and Taubes [Math. Res. Lett. 3 (1996) 661-674].}

\begin{abstract}  
Let $X$ be a closed manifold with $\chi(X)=0$, and let $f\co X\to S^1$ be
a circle-valued Morse function.  We define an invariant $I$ which
counts closed orbits of the gradient of $f$, together with flow lines
between the critical points.  We show that our invariant equals a form
of topological Reidemeister torsion defined by Turaev
\cite{turaev:spinc}.

We proved a similar result in our previous paper \cite{hutchings-lee},
but the present paper refines this by separating closed orbits and
flow lines according to their homology classes.  (Previously we only
considered their intersection numbers with a fixed level set.)  The
proof here is independent of the proof in \cite{hutchings-lee}, and
also simpler.

Aside from its Morse-theoretic interest, this work is motivated by the fact
that when $X$ is three-dimensional and $b_1(X)>0$, the invariant $I$ equals a
counting invariant $I_3(X)$ which was conjectured in \cite{hutchings-lee}
to equal the Seiberg--Witten invariant of $X$.  Our result, together with this
conjecture, implies that the Seiberg--Witten invariant equals the
Turaev torsion.  This was conjectured by Turaev \cite{turaev:spinc}
and refines the theorem of Meng and Taubes \cite{meng-taubes}.
\end{abstract}


\maketitlepage

\section{Introduction}
\label{sec:intro}

Given a flow on a manifold, it is natural to ask how many closed
orbits there are.  It turns out that for some well-behaved flows, the
numbers of closed orbits in different homology classes are related to
the Reidemeister torsion of the underlying manifold.  For example,
Fried \cite{fried:homological} defined a ``twisted Lefschetz zeta
function''
counting closed orbits of certain nonsingular hyperbolic flows and
showed that it equals a version of topological Reidemeister torsion,
which is independent of the flow.

In this paper, we are interested in the gradient flow of a
circle-valued Morse function.  For singular flows such as this one,
the zeta function is no longer invariant under deformation of the
flow.  It turns out that this lack of invariance can be fixed by considering 
the Novikov complex, which counts gradient flow lines between critical
points.  We will show that one can obtain a topological invariant by
multiplying the zeta function by the Reidemeister torsion of the
Novikov complex.  We call the resulting invariant $I$.

In our previous work \cite{hutchings-lee}, we defined a weaker version
of $I$ and showed that it equals a form of topological Reidemeister
torsion. Later we received a preprint from Turaev
\cite{turaev:spinc} defining a refined version of the latter
invariant, which we call ``Turaev torsion'' here.  Along similar lines
we can refine the Morse theoretic invariant in \cite{hutchings-lee}
to obtain the invariant $I$.
The main result of this paper asserts that $I$ equals Turaev torsion.

Our previous methods are not quite sufficient to prove
this refinement, so here we introduce a different and simpler approach.
This paper is independent of \cite{hutchings-lee}, except
that the latter paper defines certain compactifications in Morse
theory which we use here, and also provides more background and context.

We now proceed to define our invariant $I$ more precisely and state
our main theorem.  We then describe the application to
three-dimensional Seiberg--Witten theory.  In \S\ref{sec:background} we
give some background definitions, and in \S\ref{sec:proof} we prove
the main theorem.  In \S\ref{sec:sw} we give more details on the
relation to Seiberg--Witten theory.

\subsection{Statement of results}

The basic setup for this paper is as follows.  Let $X$ be a closed
connected oriented $n$--dimensional manifold.  We assume throughout
that $\chi(X)=0$, so that we can define Reidemeister torsion.  Also,
our result is most interesting when $b_1(X)>0$.

Let $f\co X\to S^1$.  In order to consider the gradient flow of $f$, we
endow $X$ with a Riemannian metric.  We make the following
assumptions:

\begin{assumption}
\label{assumption:admissible}
\begin{itemize}
\item[(a)] $f$ is a Morse function.
\item[(b)] The ascending and descending manifolds of the critical
points of $f$ intersect transversely (see \S\ref{sec:morse}).
\item[(c)] The closed orbits of the gradient $\nabla f$ are nondegenerate
(see below).
\end{itemize}
\end{assumption}
A standard transversality argument shows that these assumptions hold
if $f$ and the metric are generic.

A {\em closed orbit\/} is a nonconstant map $\gamma\co S^1\to X$ with
$\gamma'(t)=-\lambda\nabla f$ for some $\lambda>0$.  We declare two
closed orbits to be equivalent if they differ by reparametrization.
The {\em period} $p(\gamma)$ is the largest integer $p$ such that
$\gamma$ factors through a $p$--fold covering $S^1\to S^1$.  A closed
orbit is {\em nondegenerate} if $\det(1-d\phi(x))\neq 0$,
where $\phi$ is the $p^{th}$ return map at a point
$x\in\gamma(S^1)\subset X$.  If so, the {\em Lefschetz sign}
$\epsilon(\gamma)$ is the sign of this determinant.

\begin{notation}
Let $H_1:=H_1(X)$.  Let $\theta\in H^1(X;\Z)$ denote the pullback by
$f$ of the ``upward'' generator of $H^1(S^1;\Z)$.

Let $\Lambda=\Nov(H_1,-\theta)$ denote the {\em Novikov ring}
\cite{novikov,hofer-salamon}, consisting of functions $H_1\to\Z$ that
are finitely supported on the set $\{h\in H_1\mid-\theta(h)\le C\}$
for each $C\in\R$.  This ring has the obvious addition, and the
convolution product.  We denote a function $a\co H_1\to\Z$
by the (possibly infinite) formal sum $\mbox{$\sum_{h\in H_1}a(h)\cdot h$}$.
\end{notation}

\begin{definition}
\cite{fried:homological,pajitnov:zeta}\qua
We count closed orbits with the {\em zeta function}, which is a
function $H_1\to \Z$ defined, in the above notation, by
\begin{equation}
\label{eqn:zetaIntrinsic}
\zeta\eqdef \exp\left(\sum_{\gamma\in\mc{O}}
\frac{\epsilon(\gamma)}{p(\gamma)}[\gamma]\right).
\end{equation}
Here $\mc{O}$ denotes the set of closed orbits, and
$[\gamma]\eqdef\gamma_*[S^1]$ is the homology class of
$\gamma$ in $H_1(X)$.
\end{definition}

A compactness argument using Assumption~\ref{assumption:admissible},
together with the observation that $-\theta ([\gamma])>0$ for all
$\gamma$, shows that $\zeta\in\Lambda\tensor\Q$.

We remark that there is also a product formula
\cite{fried:survey,hutchings-lee,ionel-parker}
\begin{equation}
\label{eqn:productFormula}
\zeta=\prod_{\gamma\in\mc{I}}(1-(-1)^{i_-}[\gamma])^{-(-1)^{i_0}}.
\end{equation}
Here $\mc{I}$ denotes the set of irreducible (period 1) closed orbits,
and $i_-(\gamma)$ and $i_0(\gamma)$ denote the numbers of real
eigenvalues of the return map in the intervals $(-\infty,-1)$ and
$(-1,1)$, respectively.  Equation \eqref{eqn:productFormula} shows
that in fact $\zeta\in \Lambda$, ie, $\zeta$ has integer
coefficients.  A third formula for the zeta function, in terms of
fixed points of return maps, is given in equation
\eqref{eqn:zetaFixedPoints}.

We now introduce a notion of topological Reidemeister torsion
following Turaev \cite{turaev:spinc}, and an analogous notion of
Morse-theoretic torsion.  Detailed definitions are given in
\S\ref{sec:torsion}.

Let $\tilde{X}$ denote the universal (connected) abelian cover of $X$,
whose automorphism group is $H_1(X)$.  A smooth triangulation of $X$
lifts to $\tilde{X}$ and gives rise to a chain complex
$C_*(\tilde{X})$ over $\Z[H_1]$.  The Reidemeister torsion of this
complex is an element of $Q(\Z[H_1])/\pm H_1$, where $Q(R)$ denotes
the total quotient ring of $R$.  The $\pm H_1$ ambiguity arises
because the Reidemeister torsion depends on a choice of ordered basis.

Turaev \cite{turaev:euler} showed that the $H_1$ ambiguity can be
resolved by the choice of an ``Euler structure''.  The space $\Eul(X)$
of Euler structures is a natural affine space over $H_1(X)$, reviewed
in \S\ref{sec:euler}.  One can also resolve the sign ambiguity by
choosing a {\em homology orientation} of $X$, ie, an orientation $o$
of $\bigoplus_iH_i(X;\Q)$ (see \cite{turaev:knot}).  
We can then define the {\em Turaev
torsion}
\begin{equation}
\label{eqn:turaevTorsion}
\tau(X;o)\co \Eul(X)\to Q(\Z[H_1]).
\end{equation}
This is an $H_1$--equivariant map which does not depend on the
triangulation.  We write $\tau(X){:=}\pm\tau(X,o)$; this is an
$H_1$--equivariant map $\Eul(X)\to Q(\Z[H_1])/\mbox{$\pm1$}$.

\begin{example}
If $X$ is the 3--manifold obtained by zero surgery on a
knot $K\subset S^3$, then for a suitable Euler structure $\xi$,
\[
\tau(X)(\xi)=\frac{\op{Alex}(K)}{(1-t)^2}
\]
where $\op{Alex}(K)\in\Z[t,t^{-1}]/\pm1$ is the Alexander polynomial of $K$
and $t$ is a generator of $H_1(X)\simeq \Z$.
\end{example}

On the Morse theory side, the Novikov complex $\Cnov_*$ is a chain
complex over the Novikov ring $\Lambda$, whose chains are generated by
critical points of the pullback of $f$ to $\tilde{X}$, and whose
boundary operator counts gradient flow lines between critical points
(see \S\ref{sec:morse}).  We can similarly define the {\em
Morse-theoretic torsion}
\begin{equation}
\label{eqn:morseTorsion}
\tau(\Cnov_*)\co \Eul(X)\to Q(\Lambda)/\pm1.
\end{equation}

\begin{definition}
Define $I\co \Eul(X)\to Q(\Lambda)/\pm1$ to be
the product of the zeta function and the Morse-theoretic torsion:
\[
I\eqdef \zeta\cdot\tau(\Cnov_*).
\]
\end{definition}

\begin{theorem}
\label{thm:invariance}
The Morse theory invariant $I$ is independent of the metric and depends
only on the homotopy class of $f$, ie the cohomology class $\theta$.
\end{theorem}

One can prove this {\em a priori\/}; see \cite{hutchings-lee} for the
rough idea and \cite{hutchings} for the details.  Although this may
help define related invariants in other contexts, in the present
context it is easier to compute $I$ directly, which will prove
Theorem~\ref{thm:invariance} {\em a posteriori\/}.  That is what we will
do in this paper.

\begin{theorem}[Main theorem] 
\label{thm:main}
Our Morse theory invariant $I$ is equal to the topological
torsion:
\[
I=i(\tau(X))
\]
as maps $\Eul(X)\to Q(\Lambda)/\pm1$.
\end{theorem}

\noindent
Here $i\co Q(\Z[H_1])\to Q(\Lambda)$ is induced by the inclusion
$\Z[H_1]\to\Lambda$.

\begin{remarks}
(1)\qua In the extreme case when there are no critical points, $X$ is a
mapping torus and this theorem reduces to an equivariant version of
the Lefschetz fixed point theorem,
cf \cite{milnor:cyclic,fried:survey}.

(2)\qua The extreme case when
$\theta=0$, so that $f$ lifts to a real-valued Morse function, is also
essentially classical (cf \cite{milnor:whitehead}), and we treat it
in \S\ref{sec:real}.  In this case the Morse-theoretic torsion is a
topological invariant; lack of invariance and existence of closed
orbits arise simultaneously when we pass from real-valued to circle-valued
Morse theory.

(3)\qua The class $\theta\in H^1(X)$, regarded as a map $H_1(X)\to\Z$,
induces a map $\Lambda\to\Z((t))$ sending $h\mapsto t^{\theta(h)}$.
This in turn induces a partially defined map $Q(\Lambda)\to\Q((t))$.
One can deduce the main result of our previous paper
\cite{hutchings-lee} by applying this map to Theorem~\ref{thm:main}.

(4)\qua The statement of Theorem~\ref{thm:main} makes sense when $df$ is
replaced by a generic closed 1--form $\mu$ and the Novikov ring is
graded by $-[\mu]$.  It seems possible to extend
Theorem~\ref{thm:main} to this case by approximating $\mu$ by closed
1--forms in rational cohomology classes, to which
Theorem~\ref{thm:main} applies.  Another proof for closed 1--forms is
given in \cite{hutchings} by first proving
Theorem~\ref{thm:invariance} for closed 1--forms, and then using this
to reduce to the real-valued case.

(5)\qua Some previous papers, such as \cite{pajitnov:toulouse}, studied the
torsion of the Novikov complex (or the Whitehead torsion, which is
sharper but only defined ``relatively'' unless the Novikov complex is
acyclic), without considering the zeta function.  In this case one can
still obtain a topological invariant by modding out by units in the
Novikov ring with leading coefficient 1.  This is useful for
understanding the obstructions to the existence of nonvanishing closed
1--forms \cite{latour,fukaya}.  However, the extra information in the
zeta function is important for the connection with Seiberg--Witten
theory below.

(6)\qua A homology orientation of $X$ can apparently remove the sign
ambiguity in $\tau(\Cnov_*)$.  However we have not checked if
Theorem~\ref{thm:main} holds with the sign ambiguity removed this way.

\end{remarks}

\subsection{Application to Seiberg--Witten theory}

We now consider the special case when $\dim(X)=3$ and $b_1(X)>0$.  Let
$\Spinc(X)$ denote the set of spin-c structures on $X$.  Given a
homology orientation $o$, the Seiberg--Witten invariant of $X$ is a
function
\[
\SW_{X,o}\co \Spinc(X)\to\Z
\]
which counts $\R$--invariant solutions to the Seiberg--Witten equations
on $X\times \R$, modulo gauge equivalence.  (See
eg \cite{lim,meng-taubes,okonek-teleman}.)

Taubes \cite{taubes:sw=gr} has shown that the SW invariant of a
symplectic four--manifold equals a ``Gromov invariant'' counting
pseudoholomorphic curves.  In \cite{hutchings-lee} we proposed that
using similar analysis, one might be able to show that the SW
invariant of a 3--manifold is equal to a Morse theory invariant
\[
I_3\co \Spinc(X)\to\Z.
\]
The invariant $I_3$ counts certain unions of closed orbits and flow
lines of the gradient vector field of a Morse function $f\co X\to S^1$
with no index 0 or 3 critical points.  We review the definition of
$I_3$ in \S\ref{sec:sw}.

\begin{conjecture}
\label{conj:main}
\cite{hutchings-lee}\qua
The Seiberg--Witten invariant agrees with our Morse theory invariant:
\[
\SW_{X,o}=\pm I_3.
\]
\end{conjecture}

\noindent
(When $b_1(X)=1$, the SW invariant also depends on a choice of
``chamber'', and in this conjecture we use the chamber determined by
$r*df$ for $r>>0$.)

\begin{remark}
\label{remark:torus}
If $f$ has no critical points, and if we arrange for $df$ to be
harmonic, then this conjecture is a corollary of Taubes'
theorem \cite{taubes:sw=gr} applied to the symplectic four--manifold
$(X\times S^1,df\wedge ds+*_Xdf)$.  Here $s$ denotes the $S^1$
coordinate.  The idea is that for a suitable homology orientation, if
$\mathfrak{s}\in\Spinc(X)$, then
\[
\SW_{X,o}(\mathfrak{s})=\SW_{X\times S^1}(\pi^*\mathfrak{s})=\Gr_{X\times
S^1}(\pi^*\mathfrak{s})=I_3(\mathfrak{s}).
\]
The first equality expresses the fact that all solutions to the SW
equations on $X\times S^1$ are $S^1$--invariant; see
\cite{okonek-teleman} for details of this equality.  The second
equality is Taubes' theorem; here $\Gr_{X\times S^1}(\pi^*\mathfrak{s})$
counts, in the sense of \cite{taubes:counting}, pseudoholomorphic
curves in a certain $S^1$--invariant homology class in $X\times S^1$.
An energy argument shows that for a suitable almost complex structure,
every such curve is a union of closed orbits of $\nabla f$ crossed
with $S^1$.  This leads to the third equality
(cf \cite[Thm. 0.1]{ionel-parker}), using the fact that $I_3$ is a
reparametrization of the zeta function in this case.

Salamon has proved a statement equivalent to
Conjecture~\ref{conj:main} in this case using a different method
\cite{salamon}.
\end{remark}

In another direction, Turaev \cite{turaev:spinc} conjectured a
combinatorial formula for the Seiberg--Witten invariant as follows.  If
$\dim(X)=3$ and $b_1(X)>1$, then for an Euler structure $\xi$, the
torsion $\tau(X;o)(\xi)$ is actually in the group ring $\Z[H_1]$.  If
$b_1=1$, then $i(\tau(X;o)(\xi))\in\Lambda$, rather than in the
quotient ring.  Given a homology orientation $o$, one can then define
a map
\[
\begin{split}
T(X;o)\co \Eul(X)&\longrightarrow \Z,\\
\xi&\longmapsto i(\tau(X;o)(\xi))(0).
\end{split}
\]
where $(0)$ indicates evaluation on $0\in H_1$.  (This depends on the
sign of $[df]$ when $b_1=1$.)  There is also a natural isomorphism
$\imath\co \Spinc(X)\to\Eul(X)$ (\cite{turaev:euler}, see
\S\ref{sec:conclusion}).

\begin{conjecture}
\label{conj:turaev}
(Turaev \cite{turaev:spinc})\qua
The Seiberg--Witten invariant agrees with the Turaev torsion:
\[
\SW_{X,o}=T(X;o)\circ\imath\co \Spinc(X)\to\Z.
\]
\end{conjecture}

This statement is a refinement of the theorem of Meng and Taubes
\cite{meng-taubes}, which gives an ``averaged'' version of this
equivalence, in which one sums over spin-c structures that differ by
torsion elements of $H^2(X;\Z)$.

The invariant $I_3$ turns out to be a reparametrization of the more
general invariant $I$.  Thus we can apply
Theorem~\ref{thm:main} to prove:

\begin{theorem}
\label{thm:sw}
Conjecture~\ref{conj:main} is equivalent to Conjecture~\ref{conj:turaev}
(modulo signs).
\end{theorem}

The detailed proof is given in \S\ref{sec:sw}.

\begin{update}
(1)\qua Three days after the first version of this paper was
posted on the internet, a preprint by Pajitnov \cite{pajitnov:zeta}
appeared, giving a result similar to Theorem~\ref{thm:main}, using
Whitehead torsion.

(2)\qua Turaev \cite{turaev:new} has shown how to refine the methods of
Meng and Taubes to prove Conjecture~\ref{conj:turaev}, modulo signs.
Together with our results, this indirectly proves
Conjecture~\ref{conj:main}.  However, one might still desire a direct
analytic proof.  The following is a summary of the situation:

\vbox{
$$
\mbox{\scriptsize Thm. \ref{thm:main}}
$$
$$
\begin{array}{rcl} \mbox{$S^1$ Morse theory} & \mbox{\large $=$} &
\mbox{Turaev torsion} \\ \mbox{\scriptsize analytic proof?}
\raisebox{-2pt}{
\begin{picture}(15,14)
\put(0,13){\line(1,-1){13}}
\put(2,14){\line(1,-1){13}}
\end{picture}}
& &
\raisebox{-2pt}{
\begin{picture}(15,14)
\put(0,1){\line(1,1){13}}
\put(2,0){\line(1,1){13}}
\end{picture}}
\mbox{\scriptsize Meng-Taubes/Turaev}
\end{array}
$$
$$
\mbox{Seiberg--Witten}
$$
}
\end{update}

\begin{acknowledgments}
This paper would not exist were it not for Taubes' work on
Seiberg--Witten and Gromov invariants.  We thank him for sharing his
ideas generously.  We also thank R Bott, R Forman, and D Salamon
for helpful conversations.
\end{acknowledgments}

\section{Background}
\label{sec:background}

We now give some necessary background.  Section~\ref{sec:morse}
reviews the definition of the Novikov complex, which counts gradient
flow lines between critical points.  Section~\ref{sec:euler} reviews
Turaev's Euler structures, which are needed for the most refined
version of Reidemeister torsion.  Section~\ref{sec:torsion} gives
the precise definitions of the versions of Reidemeister torsion that we
use.  Finally, section~\ref{sec:real} proves Theorem~\ref{thm:main} for
real-valued Morse functions, as a warmup for some of the arguments in
\S\ref{sec:proof}.

\subsection{The Novikov complex}
\label{sec:morse}

We begin with some standard definitions from Morse theory.
If $p$ is a critical point in $X$ of $f$ or in $\tilde{X}$ of the
pullback of $f$, the {\em descending manifold} $\mc{D}(p)$ is the set
of all points $x$ such that upward gradient flow starting at $x$
converges to $p$.  Similarly the {\em ascending manifold} $\mc{A}(p)$
is the set of all points from which downward gradient flow converges
to $p$.  If $\op{ind}(p)=i$, then $\mc{D}(p)$ and $\mc{A}(p)$ are
embedded open balls of dimensions $i$ and $n-i$ respectively.

The {\em Novikov complex} $(\Cnov_*,\partial^f)$ is defined as
follows.  On $\tilde{X}$, we can lift our Morse function to a map
$\tilde{f}\co \tilde{X}\to\R$.  Let $\Cnov_i$ denote the set of (possibly
infinite) linear combinations $\alpha$ of critical points of index $i$
in $\tilde{X}$, such that for each $C\in\R$, the sum $\alpha$ contains only
finitely many critical points $p\in\tilde{X}$ with $\tilde{f}(p)>C$.
The action of $H_1$ on the critical points by covering transformations
induces an action of the Novikov ring $\Lambda$ on $\Cnov_i$.  In fact
$\Cnov_i$ is a free $\Lambda$--module; one can specify a basis for
$\Cnov_*$ by choosing a lift of each critical point in $X$ to
$\tilde{X}$.

We define $\partial^f\co \Cnov_i\to\Cnov_{i-1}$ as follows.  If
$p,q\in\tilde{X}$ are critical points of index $i$ and $i-1$
respectively, let $\langle\partial^fp,q\rangle$ denote the signed
number of gradient flow lines from $p$ to $q$.  If $p$ is a critical
point of index $i$, define
\[
\partial^fp\eqdef\sum_q\langle \partial^fp,q\rangle q
\]
where the sum is over all critical points $q\in\tilde{X}$ of index
$i-1$.  We count flow lines using the sign conventions from
\cite{hutchings-lee}.  These conventions are chosen so that
$(\partial^f)^2=0$ and so that equation \eqref{eqn:bdp=dbp} holds.

\begin{theorem}[Novikov]
\label{thm:novikov}
The homology of the Novikov complex is naturally
isomorphic to the homology of the ``half-infinite'' chains in
$\tilde{X}$:
\[
H_*(\Cnov_*,\partial^f)\simeq H_*(C_*(\tilde{X})\tensor\Lambda).
\]
\end{theorem}

See eg \cite{novikov,pajitnov:toulouse,pozniak,hutchings-lee}.  For
example, if $X=S^1$ and $f\co S^1\to S^1$ has nonzero degree, then the
homology of the Novikov complex vanishes.

\subsection{Euler structures}
\label{sec:euler}

We now discuss three different notions of ``Euler structure'' and how
they relate.  One can ignore this material at the expense of allowing
an $H_1$ ambiguity in Reidemeister torsion.

\begin{definition}
\label{def:SES1}
(Turaev \cite{turaev:euler})\qua If $X$ is a closed smooth manifold with
$\chi(X)=0$ and $n=\dim(X)>1$, a {\em smooth Euler structure} on $X$
is a nonsingular vector field on $X$, where two such vector fields are
declared equivalent if their restrictions to the complement of a ball
in $X$ are homotopic through nonsingular vector fields.
We let $\Eul(X)$ denote the space of smooth Euler structures.  
\end{definition}

By
obstruction theory, $\Eul(X)$ is an affine space over
$H^{n-1}(X;\pi_{n-1}(S^{n-1}))=H_1(X)$.  (It is nonempty since
$\chi(X)=0$.)

The following alternate definition of smooth Euler structures is
useful for Morse theory, and also works well when $n=1$.  If $v$ is a
vector field on $X$ with nondegenerate zeroes, let $H_1(X,v)$ denote
the set of homology classes of 1--chains $\gamma\in X$ with
$\partial\gamma=v^{-1}(0)$, where the points in $v^{-1}(0)$ are oriented
in the standard way.  The set $H_1(X,v)$ is a subset of the relative
homology $H_1(X,v^{-1}(0))$ and is an affine space over $H_1(X)$.

\begin{definition}
\label{def:SES2}
One can show by a cobordism argument that the spaces\break $H_1(X,v)$ for
different $v$'s are canonically isomorphic to each other, and hence to
a single space.  We call this space $\Eul(X)$.  We let
$i_v\co H_1(X,v)\to\Eul(X)$ denote the canonical isomorphism.
\end{definition}

If $n>1$, we can go from Definition~\ref{def:SES2} to
Definition~\ref{def:SES1} as follows. Given $\gamma\in H_1(X,v)$, we
can represent $\gamma$ by disjoint paths connecting the zeroes of $v$
in pairs.  We then construct a nonsingular vector field by cancelling
the zeroes of $v$ in a neighborhood of $\gamma$.  (If $v$ has no
zeroes, we send $0\in H_1(X)=H_1(X,v)$ to the Euler structure
represented by $v$ and extend equivariantly.)

\begin{definition}\label{def:ceuler}
\cite{turaev:euler}\qua
Let $(X,T)$ be a finite connected CW--complex with cells $\{\sigma_i\}$.
($X$ denotes the underlying topological space; $T$ denotes the cell structure.)
A {\em combinatorial Euler structure} on $(X,T)$ is a
choice of a lift of each cell to the universal abelian cover
$\tilde{X}$, where two such sets of lifts
$\{\tilde{\sigma_i}\}$ and $\{h_i\tilde{\sigma_i}\}$, with $h_i\in
H_1(X)=\op{Aut}(\tilde{X})$, are considered equivalent if
$\sum_i(-1)^{\dim(\sigma_i)}h_i=0$.
\end{definition}

We let $\Eul(X,T)$ denote the
space of combinatorial Euler structures of the CW--complex $(X,T)$.
This is clearly an affine space over $H_1(X)$.

Note that if $\bar{T}$ is a refinement of the cell-structure $T$ with
cells $\{\tau_j\}$, then there is a canonical isomorphism
$\Eul(X,T)\to \Eul(X,\bar{T})$ sending $\{\tilde{\sigma}_i\}$ to
$\{\tilde{\tau}_j\}$, where $\tilde{\tau}_j\subset\tilde{\sigma}_i$ if
$\tau_j\subset \sigma_i$.

\begin{lemma} {\rm\cite{turaev:euler}}\qua
\label{lem:euler}
If $X$ is a closed smooth manifold with a smooth triangulation $T$,
then there is a natural isomorphism between the spaces of smooth and combinatorial
Euler structures:
\[
\Eul(X)\simeq\Eul(X,T).
\]
\end{lemma}

The idea is that there is a natural vector field on each simplex
with a zero at the center of each face and which points into the
simplex near the boundary.  These piece together to give a continuous
vector field on $X$.  We can perturb this to a smooth vector field
$v_T$ with a nondegenerate zero of sign $(-1)^i$ in the center of each
$i$--simplex.  Then a smooth Euler structure $\xi$ can be represented
by a chain $\gamma$ consisting of paths connecting the zeroes in
pairs, with $[\gamma]=i_{v_T}^{-1}(\xi)\in H_1(X,v_T)$.  We can lift
the chain $\gamma$ to $\tilde{X}$, and the induced lifts of its
endpoints determine a combinatorial Euler structure.

\subsection{Reidemeister torsion}
\label{sec:torsion}

We now review the definition of Reidemeister torsion of certain chain
complexes.  We then use this algebra to define Reidemeister torsion
for the two geometric complexes we are interested in.

\subsubsection{Algebra}
\label{sec:torsionAlgebra}

Let $(C_*,\partial)$ be a finite complex of finite dimensional
vector spaces over a field $F$.  The standard short exact sequences
$0\to Z_i\to C_i\to B_{i-1}\to 0$ and $0\to B_i\to Z_i\to H_i\to 0$
induce canonical isomorphisms
$\det(C_i)=\det(Z_i)\tensor\det(B_{i-1})$ and
$\det(Z_i)=\det(B_i)\tensor\det(H_i)$, where `$\det$' denotes top
exterior power.  Combining these isomorphisms gives an isomorphism
\begin{equation}
\label{eqn:determinantIsomorphism}
\Phi\co\bigotimes_i\det(C_i)^{(-1)^i}\to \bigotimes_i\det(H_i)^{(-1)^i}.
\end{equation}

Let $e$ be an ordered basis for $C_*$, ie, an ordered basis $e_i$ for
each $C_i$.  Let $h$ be an ordered basis for $H_*$.
Let $[e]\in \bigotimes_i\det(C_i)^{(-1)^i}$ and
$[h]\in\bigotimes_i\det(H_i)^{(-1)^i}$ denote the resulting volume forms.

\begin{definition}
We define the {\em Reidemeister torsion}
\[
\hat{\tau}(C_*,e,h)\eqdef\Phi([e])/[h]\in F^\times.
\]
We also define
\[
\tau(C_*,e)\eqdef\left\{\begin{array}{cl}\hat{\tau}(C_*,e,1) & \mbox{if
$H_*=0$,}\\ 0 & \mbox{otherwise.}\end{array}\right.
\]
\end{definition}

Usually we will be interested in $\tau$ rather than $\hat{\tau}$.  In
practice, we can compute the torsion $\tau$ in terms of an alternating
product of determinants as follows.

\begin{lemma}
\label{lem:computeTorsion}
If $H_*=0$, we can find a
decomposition $C_*=A_*\oplus B_*$ such that (i) $A_i$ and $B_i$ are
spanned by subbases of the basis $e_i$, and (ii) the map
\[
\pi_{B_{i-1}}\circ\partial|_{A_i}\co A_i\to B_{i-1}
\]
(which we abbreviate by $\partial\co A_i\to B_{i-1}$) is an isomorphism.
Then
\[
\tau(C_*,e)\eqdef\pm\prod_i\det(\partial\co A_i\to B_{i-1})^{(-1)^i}.
\]
{\em
Here the determinants are computed using the subbases of $e$.}
\end{lemma}

We now extend the definition of torsion to complexes over
certain rings which might not be fields.

\begin{definition}
\cite{turaev:spinc}\qua
\label{def:rings}
Let $R$ be a ring, and assume that its total quotient ring (denoted by
$Q(R)$) is a
finite sum of fields, $Q(R)=\oplus_jF_j$.  Let $(C_*,\partial)$ be a
finite complex of finitely generated free $R$--modules with
an ordered basis $e$. We define
\[
\tau(C_*,e)\eqdef\sum_j\tau(C_*\tensor_RF_j,e\tensor 1)\in\bigoplus_jF_j=Q(R).
\]
In this case $H_*(C)$ might not be free, in which case it does not
have a basis in the usual sense. However in this paper we call a set
$h:=\{h_j\}$ a ``basis'' for $H_*(C)$ when $h_j$ is a basis for
$H_*(C_*\tensor_RF_j)$ for each $j$.  Given $h=\{h_j\}$, we define
\[
\hat{\tau}(C_*,e,h)\eqdef\sum_j \hat{\tau}(C_*\tensor_RF_j,e\tensor
1,h_j)\in Q(R)^\times.
\]
\end{definition}

\begin{example}
\label{example:determinant}
If $0\to
C_2\stackrel{\partial}{\to}C_1\to 0$ is a 2--term complex with
$C_1= C_2$, and if $e$ is a basis which is identical on $C_1$ and
$C_2$, then
\[
\tau(C_*,e)=\det(\partial).
\]
\end{example}

We are interested in the rings $\Z[H_1]$ and $\Lambda$. Their
quotients are finite sums of fields (see eg \cite{turaev:spinc}), and
these decompositions are compatible with the inclusion
$\Z[H_1]\to\Lambda$.

In \S\ref{sec:technical}, we will need the following product formula
for torsion.  Let $R$ be a ring such that $Q(R)$ is a finite sum of
fields $F_j$.  Let $0\to C_*\to C_*'\to C_*''\to 0$ be a short exact
sequence of finite complexes of finitely generated free $R$--modules.
Let $e,e',e''$ be bases for $C_*,C_*',C_*''$ compatible with the exact
sequence.  Let $h,h',h''$ be bases for the homology as in
Definition~\ref{def:rings}.  Let $L_*$ denote the long exact sequence
in homology, regarded as an acyclic chain complex, and let $b$ denote
the basis for $L_*$ obtained by combining $h,h',h''$.

\begin{lemma}
\label{lem:productFormula}
We have the following product formula for torsion:
\[
\hat{\tau}(C',e',h')=\hat{\tau}(C,e,h)\hat{\tau}(C'',e'',h'')\tau(L_*,b).
\]
\end{lemma}

\begin{proof}
This follows from \cite{milnor:whitehead}.
\end{proof}

\subsubsection{Geometry}
\label{sec:torsionGeometry}

We now define the Turaev torsion \eqref{eqn:turaevTorsion}.  In the
notation of Definition \ref{def:ceuler}, let $(X,T)$ be a finite
connected CW--complex with universal abelian cover $\tilde{X}$.  Lifing
the cells gives a chain complex $C_*(\tilde{X},T)$ over $\Z[H_1(X)]$.
A combinatorial Euler structure $\xi$ determines a set of lifts of
each cell to $\tilde{X}$, up to equivalence. Choose one set of lifts;
this gives a basis for $C_*(\tilde{X},T)$.  A homology orientation $o$
determines an orientation of this basis, via the isomorphism
\eqref{eqn:determinantIsomorphism} applied to $C_*(X,T)$.  Let
$e(\xi,o)$ denote the resulting ordered basis.  We define the {\em
combinatorial Turaev torsion} $\tau(X,T;o)$ to be the
$H_1$--equivariant map $\Eul(X,T)\to Q(\Z[H_1])$ given by
\[
\tau(X,T;o)(\xi)\eqdef\tau(C_*(\tilde{X},T),e(\xi,o))\in Q(\Z[H_1]).
\] 
Note that the right hand side of this equation does not depend on the choice of
a set of lifts. Furthermore, $\tau(X,T;o)=\tau(X,\bar{T};o)$ under the 
canonical isomorphism $\Eul(X,T)\to\Eul(X, \bar{T})$, if $\bar{T}$ is a 
refinement of $T$.

\begin{definition}
Let $X$ be a closed connected smooth manifold with $\chi(X)=0$, with a
smooth Euler structure
$\xi\in\Eul(X)$ and a homology orientation $o$.  Choose a smooth
triangulation $T$ of $X$.  Let $\xi_T\in\Eul(X,T)$ denote the combinatorial
Euler structure equivalent to $\xi$ via Lemma~\ref{lem:euler}.  We
define the {\em Turaev torsion}
\[
\tau(X;o)(\xi)\eqdef\tau(X,T;o)(\xi_T).
\]
\end{definition}

The results of \cite{turaev:euler} show that the Turaev torsion does
not depend on the choice of smooth triangulation $T$.

We now define the Morse theoretic torsion
\eqref{eqn:morseTorsion}.  A smooth Euler structure $\xi$ can be represented
by a chain $\gamma$ connecting the critical points of $f$ in pairs, with
$[\gamma]=i_{\nabla f}^{-1}(\xi)\in H_1(X,\nabla f)$.  We can lift $\gamma$ to
$\tilde{X}$, and the induced lifts of the endpoints determine a basis
$e(\xi)$ for $\Cnov_*$.

\begin{definition}
We define the {\em Morse theoretic torsion} $\tau (CN_*)\co  \mbox{Eul}(X)\to Q(\Lambda)/\pm 1$ by
\[
\tau(CN_*)(\xi)\eqdef\tau(CN_*,e(\xi)).
\]
\end{definition}
The map $\tau(CN_*)$ is $H_1$--equivariant, and again does not depend
on the choice of lifting.  There is a sign ambiguity because the basis
$e(\xi)$ is unordered. In the special case when $f$ has no critical
points, we define $\tau(\Cnov_*)$ to be the $H_1$ equivariant map such
that $\tau(\Cnov_*)(\xi)=1$ for the smooth Euler structure $\xi$
represented by $\nabla f$.

In the future, we call two bases for $C_*(\tilde{X},T)$ or $\Cnov_*$ equivalent
if they correspond to the same Euler structure $\xi$.

\subsection{The real-valued case}
\label{sec:real}

Before proceeding more deeply into circle-valued Morse theory, it will be
useful to prove the main theorem for real-valued Morse functions.

\begin{lemma}
\label{lem:real}
Theorem~\ref{thm:main} holds when
$f\co X\to S^1$ lifts to a real-valued Morse function $X\to\R$.
\end{lemma}

\begin{proof}
In this case $\zeta=1$, so we just need to check that the Morse
theoretic and topological torsions agree.  This is essentially
classical (cf \cite{milnor:whitehead}), except for the identification
of the bases determined by an Euler structure.

If $\xi$ is an Euler structure, then the bifurcation analysis in
\cite{laudenbach} shows that $\tau(CN_*)(\xi)$ is independent of the
real-valued Morse function and the metric.  It is not hard to check
that the Euler structures work out at each stage.  In \cite{laudenbach}
it is assumed that the metric has a standard form near the critical
points, but this can be arranged by a perturbation which does not
affect the Novikov complex.

Now let $T$ be a smooth triangulation of $X$.  We can apparently
perturb the vector field $V_T$ of \S\ref{sec:euler} so that it is the
gradient of a Morse function $F\co X\to\R$ with respect to some metric.
In this case, the Novikov complex $(\Cnov_*,\partial^{F})$ is
identical to the chain complex $(C_*(\tilde{X};T),\partial)$.
Moreover the bases determined by $\xi$ agree.  Thus the
Morse-theoretic and topological torsion are equal.
\end{proof}

\section{Proof of the main theorem}
\label{sec:proof}

We will now prove the main theorem as follows.  In \S\ref{sec:complex}
we prepare for the computation of torsion by constructing a cell
complex $X'$ which ``approximates'' $X$ and is adapted to the vector
field $\nabla f$.  In \S\ref{sec:technical} we prove a technical lemma
(Lemma~\ref{lem:technical}) asserting that $X$ and $X'$ have the same
Reidemeister torsion.  The heart of the proof is in
\S\ref{sec:compute} and \S\ref{sec:interpret}, where we determine the
torsion of $X'$ by a short computation, and then interpret the answer
geometrically to recover the invariant $I$.

\subsection{The cell complex $X'$}
\label{sec:complex}

Assume $0\in S^1=\R/\Z$ is a regular value of $f$ (by composing $f$
with a rotation if necessary).  Let $\Sigma\eqdef f^{-1}(0)$.  Let $Y$
be the compact manifold with boundary obtained by cutting $X$ along
$\Sigma$.  We can write $\partial Y=\Sigma_1\sqcup\Sigma_0$, where
$\Sigma_i$ is canonically isomorphic to $\Sigma$, and $-\nabla f$
points inward along $\Sigma_1$.

We give $Y$ a cell decomposition as follows.  Let $T_1$ be a smooth
triangulation of $\Sigma_1$ such that each simplex is transverse to
the ascending manifolds of the critical points in $Y$.  If $p\in Y$ is
a critical point, let $\mc{D}_0(p)$ denote the descending manifold of
$p$ in $Y$.  If $\sigma\in T_1$ is a simplex, let $\mc{F}(\sigma)$
denote the set of all $y\in Y$ such that upward gradient flow starting
at $y$ hits $\sigma$.  Choose a cell decomposition $T_0$ of
$\Sigma_0$, such that the intersections with $\Sigma_0$ of
$\mc{D}_0(p)$ and $\mc{F}(\sigma)$ are subcomplexes, for each critical
point $p$ and each simplex $\sigma\in T_1$.

\begin{lemma}
\label{lem:cellY}
The cells in $T_1$ and $T_0$, together with
all the cells $\mc{D}_0(p)$ and $\mc{F}(\sigma)$, give a legitimate cell
decomposition, $T'_Y$, of $Y$.
\end{lemma}

\begin{proof}
Recall that
$\mc{D}_0(p)$ and $\mc{F}(\sigma)$ have natural compactifications
using broken flow lines (cf \cite{hutchings-lee}). It may be shown
by ``induction on height'' that these compactifications are
homeomorphic to closed balls.  There are moreover natural continuous
maps of the compactifications to $Y$ which send the interiors of the
balls homeomorphically to $\mc{D}_0(p)$ and $\mc{F}(\sigma)$.  The
transversality condition on $T_1$ and Assumption~\ref{assumption:admissible}(b) ensure that the boundary of a cell consists
of lower dimensional cells in $T'_Y$.
\end{proof}

We would like to glue the boundary components of $Y$ back together to
obtain a nice cell decomposition of $X$,
but usually $T_0$ will not agree with $T_1$.  To correct for this, let
$\rho\co (\Sigma_0,T_0)\to(\Sigma_1,T_1)$ be a cellular approximation to
the canonical identification $\Sigma_0\to\Sigma_1$.  Consider the
mapping cylinder of $\rho$:
\[
M_\rho=\frac{(\Sigma_0\times[0,1])\sqcup \Sigma_1}{(x,1)\sim \rho(x)}.
\]
This has a cell decomposition consisting of $T_0$ and $T_1$, together
with the cells $\Delta\times(0,1)$ for each $\Delta\in T_0$.
There is a canonical inclusion $\Sigma_0\to M_\rho$ sending $x\mapsto
(x,0)$, and there is also a canonical inclusion $\Sigma_1\to M_\rho$.

\begin{definition}
Let $X'$ be the space obtained by gluing $Y$ and $M_\rho$ along
$\Sigma_0\sqcup\Sigma_1$.
\end{definition}

The space $X'$ inherits a cell decomposition, but for our computations
we prefer a simpler cell decomposition, obtained by fusing some cells
together as follows.  If $\Delta$ is a cell in $Y$ of the form
$\mc{D}_0(p)$ or $\mc{F}(\sigma)$, we define a corresponding cell in $X'$ by
\[
\hat{\Delta}\eqdef\Delta\cup ((\partial\Delta\cap \Sigma_0)\times[0,1)).
\]
Here $(\partial\Delta\cap \Sigma_0)\times[0,1)$
indicates a subset of $M_{\rho}$.

\begin{definition}
\label{def:fusion}
Let $T'$ be the cell decomposition of $X'$ consisting of cells of the
following types:
\begin{itemize}
\item[(a)]
$\widehat{\mc{D}_0(p)}$ for $p\in Y$ a critical point;
\item[(b)]
simplices in $T_1$;
\item[(c)]
$\widehat{\mc{F}(\sigma)}$ for $\sigma\in T_1$.
\end{itemize}
\end{definition}

\subsection{$X$ and $X'$ have the same Reidemeister torsion}
\label{sec:technical}

We now show that $X$ and $X'$ have the same Reidemeister torsion, if
the Euler structures are compatible in an appropriate sense.

We begin by noting that $H_*(X')=H_*(X)$, and 
$H_*(\tilde{X'})=H_*(\tilde{X})$ as $\Z[H_1]$--modules, 
as one can see from the
exact sequences \eqref{eqn:exactsq-x}, \eqref{eqn:exactsq-x'} below.
Note that the universal abelian cover
$\tilde{X}'$ of $X'$ is obtained from $\tilde{X}$ by modifying a
neighborhood of the inverse image of $\Sigma$. 

\begin{notation}
(1)\qua If $Z$ is a subset of $X$ or $X'$, then $\tilde{Z}$ will denote the
inverse image of $Z$ in $\tilde{X}$ or $\tilde{X'}$.  So $\tilde{Z}$
is usually not the universal abelian cover of $Z$.

(2)\qua  We omit the cell structures from the notation when they are clear
    from context.
\end{notation}

A smooth Euler structure $\xi$ on $X$ corresponds to an equivalence
class of lifts of the critical points of $f$ to $\tilde{X}$, as in
\S\ref{sec:torsionGeometry}.  A combinatorial Euler structure $\xi'$ on $X'$
consists of an equivalence class of lifts of the cells $T'$ to
$\tilde{X'}$.

\begin{definition}
\label{def:compatible}
We say that $\xi$ and $\xi'$ are {\em compatible} if,
within these equivalence classes, the lifts can be chosen so that:
\begin{itemize}
\item[(a)] The lift of each critical point $p$ in $\xi$ is contained in
the lift of the cell $\widehat{\mc{D}_0(p)}$ in $\xi'$.
\item[(b)]
For each simplex $\sigma\in T_1$, the lift of $\sigma$ in $\xi'$
agrees with the ``top'' of the lift of $\widehat{\mc{F}(\sigma)}$ in $\xi'$.
\end{itemize}
\end{definition}
The compatibility conditions in Definition
\ref{def:compatible} induce an isomorphism from $\mbox{Eul}(X)$ to
$\mbox{Eul}(X',T')$ as affine spaces over $H_1$.

Recall that $\tau(X):=\pm \tau(X,o)$; similarly write $\tau(X',T'):=\pm
\tau(X',T',o)$.

\begin{lemma}
\label{lem:technical}
If the Euler structures $\xi\in \mbox{Eul}(X)$ and $\xi'\in
\mbox{Eul}(X',T')$ are compatible as above, then
\[
\tau(X)(\xi)=\pm\tau(X',T')(\xi').
\]
\end{lemma}

\begin{proof}
The strategy is to compute the torsion of $X$ and $X'$ by cutting them
into pieces and using the product formula
(Lemma~\ref{lem:productFormula}) applied to various exact sequences,
and see that we obtain the same answer.  We proceed in three steps.

{\bf Step 1}\qua Consider the cell decomposition on $\tilde{\Sigma}\times
[0,1]$ consisting of cells $\tilde{\Delta}\times\{0\}$,
$\tilde{\Delta}\times(0,1)$, and $\tilde{\Delta}\times\{1\}$, where
$\tilde{\Delta}$ is a lift of a simplex $\Delta\in T_1$.  Also recall that
$\tilde{M_\rho}$ has a natural cell decomposition. 
We claim that with respect to these cell structures,
\begin{equation}
\label{eqn:step1}
\hat{\tau}(C_*(\tilde{\Sigma}\times [0,1]),e_{\Sigma},h) =
\pm\hat{\tau}(C_*(\tilde{M_\rho}),e'_{\Sigma},h'),
\end{equation}
provided that the bases $e_{\Sigma},h,e'_{\Sigma},h'$ satisfy the following
conditions:
\begin{itemize}
\item[(a)] The bases $h,h'$ for homology agree under the isomorphism induced
by the canonical map $\Sigma\times [0,1]\to M_\rho$.
\item[(b)]
The bases $e_{\Sigma},e'_{\Sigma}$ are given by lifts of cells such that:
\begin{itemize}
\item[(i)]
The lifts of the cells in
$\Sigma\times\{1\}\subset\Sigma\times[0,1]$ and $\Sigma_1\subset
M_\rho$ agree.
\item[(ii)] The lift of each cell $\Delta$ in
$\Sigma\times\{0\}\subset\Sigma\times[0,1]$ or in
$\Sigma_0\times\{0\}\subset M_\rho$ is adjacent to the lift of the cell
$\Delta\times(0,1)$ in $\Sigma\times[0,1]$ or $M_{\rho}$ respectively.
\end{itemize}
\end{itemize}

To prove \eqref{eqn:step1}, we compute both sides by applying the product
formula for torsion to the relative exact sequences
\[
0\to C_*(\tilde{\Sigma}\times\{1\})\to
C_*(\tilde{\Sigma}\times[0,1])\to
C_*(\tilde{\Sigma}\times[0,1],\tilde{\Sigma}\times\{1\})\to 0,
\]
\[
0\to C_*(\tilde{\Sigma_1})\to C_*(\tilde{M_\rho})\to
C_*(\tilde{M_\rho},\tilde{\Sigma_1})\to 0.
\]
The answers agree, because condition (a) implies that the $\tau(L_*)$
factors agree, condition (b(i)) implies that the
$\hat{\tau}(C_*(\tilde{\Sigma}))$ factors agree, and condition (b(ii))
implies that
\[
\hat{\tau}(C_*(\tilde{\Sigma}\times[0,1],\tilde{\Sigma}\times\{1\}),e_\Sigma,1)
=\pm\hat{\tau}(C_*(\tilde{M_\rho},\tilde{\Sigma_1}),e'_\Sigma,1)=\pm 1.
\]

{\bf Step 2}\qua  Let $T_Y$ be a smooth triangulation on $Y$ whose restriction
to each component of $\partial Y=\Sigma_0\sqcup\Sigma_1$ agrees with $T_1$.  
The smooth Euler structure $\xi$ on $X$ determines an equivalence class of 
bases, $e_Y$, for $C_*(\tilde{Y},\tilde{\Sigma}_0;T_Y)$, because $T_Y$ glues 
to a smooth triangulation of $X$.

Let $T_Y'$ denote the cell decomposition of $Y$
given by Lemma~\ref{lem:cellY}.  The combinatorial Euler structure $\xi'$ on $X'$ determines an equivalence class of bases $e'_Y$ for 
$C_*(\tilde{Y},\tilde{\Sigma_0};T_Y')$, because the cells of $T_Y'$ in 
$Y\setminus\Sigma_0$ are in one to one
correspondence with the cells of $X'$.

We claim that
\begin{equation}
\label{eqn:step2}
\hat{\tau}(C_*(\tilde{Y},\tilde{\Sigma_0};T_Y'),e'_Y,h')
=\pm\hat{\tau}(C_*(\tilde{Y},\tilde{\Sigma_0};T_Y),e_Y,h)
\end{equation}
provided that the bases $h'$ and $h$ on homology agree.

To prove \eqref{eqn:step2}, note that the pullback of the Morse
function $f$ to $Y$ lifts to a real-valued function $\hat{f}\co Y\to\R$.
Let $CM_*(\hat{f})$ denote the Morse complex of $\hat{f}$ on the
covering $\tilde{Y}$. A direct computation, using the compatibility of $\xi$ and $\xi'$, shows that
\[
\hat{\tau}(C_*(\tilde{Y},\tilde{\Sigma_0};T_Y'),e'_Y,h')
= \pm \hat{\tau}(CM_*(\hat{f}),e'',h'').
\]
(For similar calculations see
\S\ref{sec:compute} and \S\ref{sec:interpret}; the result here corresponds
essentially to setting $t=0$ in \eqref{eqn:bigMatrix} and Lemma~\ref{lem:compute}.) 
Here the basis $e''$ for $CM_*(\hat{f})$ is determined by the lifts of the
critical points determined by $\xi$ as before, and we assume that the
bases $h',h''$ on homology agree under the standard isomorphism
$H_*(\tilde{Y},\tilde{\Sigma_0})\simeq H_*(CM_*(\hat{f}))$.

We also have
\[
\hat{\tau}(CM_*(\hat{f}),e'',h'')
 = \pm \hat{\tau}(C_*(\tilde{Y},\tilde{\Sigma_0};T_Y),e_Y,h).
\]
The idea of the proof is to vary $\hat{f}$ in the space of Morse
functions on $Y$ such that (i) the gradient points outward along $\Sigma_0$
and does not point outward along $\Sigma_1$, and (ii) wherever the
gradient is tangent to $\Sigma_1$, the inward covariant derivative of
the gradient points inward. As in
\S\ref{sec:real}, one can show that the resulting torsion is
independent of the Morse function. Deforming $\hat{f}$ to a
Morse function $F_Y$ adapted to the triangulation $T_Y$, such that 
$CM_*(F_Y)=C_*(\tilde{Y}, \tilde{\Sigma}_0,T_Y)$, we have
$\hat{\tau}(CM_*(\hat{f}),e'',h'')=\pm \hat{\tau}(CM_*(F_Y),i_{fF}(e''),h'')$,
where $i_{fF}(e'')$ is the equivalence class of lifts (which correspond to
bases of $CM_*$) induced from $e''$ via the homotopy from $\hat{f}$ to $F_Y$.  
But $i_{fF}(e'')=e_Y$ because the homotopy from $\nabla\hat{f}$ to $\nabla F_Y$
extends to a homotopy of vector fields on $X=Y\cup (\Sigma\times I)$, which is homotopic to a homotopy from $\nabla f$ to the standard vector field associated to the triangulation of $X$
(cf end of \S2.2).

The above two equations prove \eqref{eqn:step2}.

{\bf Step 3}\qua  We now use \eqref{eqn:step1} and \eqref{eqn:step2} to compute the torsion of $X$ and $X'$.

We can regard $X$ as the union of $Y$ and
$\Sigma\times[0,1]$ along $\Sigma_0\sqcup\Sigma_1$.  Let
$\overline{T}$ denote the cell decomposition of $X$ obtained by gluing
the triangulation $T_Y$ of $Y$ to the product cell structure on
$\Sigma\times[0,1]$ obtained from $T_1$.  We then have a short exact
sequence
\begin{equation}
\label{eqn:exactsq-x}
0\to C_*\left(\tilde{\Sigma_0}\sqcup\tilde{\Sigma_1}\right)\to
C_*(\tilde{\Sigma}\times[0,1])\oplus C_*(\tilde{Y};T_Y)\to
C_*(\tilde{X};\overline{T})\to 0.
\end{equation}

Let $\overline{T'}$ denote
the ``unfused'' cell decomposition of $X'$ from \S\ref{sec:complex}.
We then have a short exact sequence
\begin{equation}
\label{eqn:exactsq-x'}
0\to C_*\left(\tilde{\Sigma_0}\sqcup\tilde{\Sigma_1}\right)\to
C_*(\tilde{M_\rho})\oplus C_*(\tilde{Y};T_Y')\to
C_*(\tilde{X'};\overline{T'})\to 0.
\end{equation}

We can choose representatives $e_Y, e'_Y$ such that they agree on $\Sigma\times
\{1\}$ and $\Sigma_1$ with $e_{\Sigma}, e'_{\Sigma}$ respectively.
Let $e(\xi)$ denote the basis for $C_*(\tilde{X};\overline{T})$
obtained by combining the bases $e_{\Sigma},e_Y$ of Step 1 and Step 2 
respectively. Similarly let
$e'(\xi')$ denote the basis for $C_*(\tilde{X'};\overline{T'})$
obtained by combining the bases $e'_{\Sigma},e'_Y$ of Step 1 and Step 2.  
Then $e(\xi)$ is a representative of the combinatorial Euler structure
on $(X,\bar{T})$ corresponding to the smooth Euler structure $\xi$, and
$e'(\xi')$ is a representative of the image of $\xi'$ under the canonical 
isomorphism $\Eul(X',T')\to \Eul(X', \bar{T}')$.

Applying the product formula to the above exact sequences, and using
equations \eqref{eqn:step1} and \eqref{eqn:step2}, we obtain
\begin{equation}
\label{eqn:almostDone}
\hat{\tau}(C_*(\tilde{X};\overline{T}),e(\xi),h)
=\pm\hat{\tau}(C_*(\tilde{X'};\overline{T'}),e'(\xi'),h').
\end{equation}
Here we are assuming that the bases $h,h'$ for homology agree under
the natural isomorphism $H_*(\tilde{X})\simeq H_*(\tilde{X'})$. 
Also, to apply \eqref{eqn:step2} in the above computation, one
relates $\tilde{Y}$ to the pair $(\tilde{Y},\Sigma_0)$ as in Step 1.

In particular, equation \eqref{eqn:almostDone} implies that
\[
\tau(C_*(\tilde{X};\overline{T}),e(\xi))
=\pm\tau(C_*(\tilde{X'};\overline{T'}),e'(\xi')).
\]
This implies the lemma because
$\tau(C_*(\tilde{X};\overline{T}),e(\xi))=\tau(X)(\xi)$, since the
insertion of $\Sigma\times[0,1]$ changes nothing, and similarly
$\tau(C_*(\tilde{X'};\overline{T'}),e'(\xi'))=\tau(X',T')(\xi')$.
\end{proof}

\subsection{Computing the torsion}
\label{sec:compute}

We now compute the torsion of the approximating cell complex $(X',T')$,
for a combinatorial Euler structure $\xi'$ compatible with a smooth
Euler structure $\xi$ on $X$ as in the previous section.

Since Lemma~\ref{lem:real} proves the theorem for real-valued Morse
functions, we assume from now on that
\[
\theta\neq 0.
\]
Without loss of generality, we may also assume that $\theta$ is
indivisible in $H^1(X;\Z)$.  (If $\theta$ is divisible by $k$, we can
lift $f$ to a $k$--fold cover of $S^1$ without changing the invariant
$I$.)  Let $V\eqdef\Ker(\theta)$, and choose a splitting
\[
H_1(X)=V\oplus \Z.
\]
Let $t$ denote the generator of the $\Z$ component with
$\theta(t)=-1$.  Then the Novikov ring can be identified with the ring
of formal Laurent series in $t$ with coefficients in $\Z[V]$:
\[
\Lambda=\Z[V]((t)).
\]

Recall that $Q(\Lambda)$ is a finite sum of fields. To prove 
Theorem~\ref{thm:main}, it suffices to show that it holds after projecting
to each such field. 
Let $K$ be a field component of $Q(\Lambda)$. By the
Novikov isomorphism (Theorem~\ref{thm:novikov}), the complexes
$\Cnov_*\tensor K$ and $C_*(\tilde{X})\tensor K$ have isomorphic
homology.  So we will assume that these complexes are both acyclic,
since otherwise they both have zero torsion $\tau$,
and there is nothing to
prove.  In all of the calculations below, we implicitly tensor everything with
the field $K$.

We can decompose
\[
C_i(\tilde{X'})=D_i\oplus E_i\oplus F_i
\]
where the three summands are generated by the cells of types (a), (b), and (c)
respectively from Definition~\ref{def:fusion}.  Let us choose a basis $e(\xi')$
for $C_i(\tilde{X'})$ as in Definition~\ref{def:compatible}. We
can identify
\[
F_i\simeq E_{i-1}.
\]
The matrix for the boundary operator on $C_i(\tilde{X'})$ can then be
written as
\begin{equation}
\label{eqn:bigMatrix}
\partial_i= \bordermatrix{ & D_i & E_i & F_i \cr
D_{i-1} &{\bf N}_i & 0 &{\bf W}_i \cr
E_{i-1} & -t{\bf M}_i & \partial^\Sigma_i & 1-t\phi_{i-1} \cr
F_{i-1} & 0 & 0 & -\partial_{i-1}^\Sigma
\cr }.
\end{equation}
We remark that $\partial_i^\Sigma$ is the boundary operator on
$C_*(\tilde{\Sigma})$.  Also $\phi_{i-1}$ is a matrix with entries in
$\Z[V]$, which can be interpreted as the return map of the gradient
flow from $\tilde{\Sigma}$ to $\tilde{\Sigma}$, after perturbation by
our cellular approximation $\rho$.  Likewise ${\bf M}_i$ sends 
$\widehat{\mc{D}_0(p)}$,
where $p\in\tilde{X}$ is a critical point, to a perturbation of
the ``descending slice'' $\widehat{\mc{D}_0(p)}\cap\tilde{\Sigma}$.

Continuing the calculation, due to the acyclicity assumption we may choose 
decompositions $D_i=D_i^A\oplus D_i^B$ such that $D_i^A$ and $D_i^B$ 
are spanned by (cells corresponding to) critical
points, and the differential $\partial^f$ induces an isomorphism
$D_i^A\to D_{i-1}^B$.  In the notation below, we denote matrices with domain
or range $D_*$ by boldface letters, and we denote their restrictions to $D_i^A$
and/or projections to $D_{i-1}^B$ by plain letters.

We now apply
Lemma~\ref{lem:computeTorsion} with $A_i=D_i^A\oplus F_i$ and
$B_i=D_i^B\oplus E_i$.  (We will explain in a moment why this choice of $A_i$ and $B_i$ is legitimate.)  We obtain
\[
\tau(C_*(\tilde{X'},e(\xi'))
=\prod_{i=1}^n\det(\Omega_i)^{(-1)^{i+1}}
\]
where
\[
\Omega_i=\bordermatrix{& D_i^A & F_i \cr
D_{i-1}^B & N_i & W_i \cr
E_{i-1} & -tM_i & 1-t\phi_{i-1} \cr}.
\]
We note that $1-t\phi_{i-1}$ is invertible because $\phi_{i-1}$ has
entries in $\Z[V]$.  It follows that
\[
\det(\Omega_i)=\det(1-t\phi_{i-1})\det(K_i)
\]
where
\begin{equation}
\label{eqn:bfK}
{\bf K}_i\eqdef {\bf N}_i+t{\bf W}_i(1-t\phi_{i-1})^{-1}{\bf
M}_i\co D_i\to D_{i-1}.
\end{equation}

It will follow from 
Lemma~\ref{lem:interpret}(b) and the choice of $D_i^A,D_i^B$ that $K_i$ is nonsingular,
provided that the triangulation $T_1$ is sufficiently fine and the
cellular approximation $\rho$ is sufficiently close to the identity.  In particular, the matrices $\Omega_i$ are then nonsingular, so that Lemma~\ref{lem:computeTorsion} legitimately applies to the $A_i$ and $B_i$ chosen above.

In conclusion, the above calculations imply the following lemma.

\begin{lemma}
\label{lem:compute}
If $T_1$ is sufficiently fine and $\rho$ is sufficiently close to the identity, then
\[
\tau(C_*(\tilde{X'}))(\xi')=
\prod_{i=1}^n\Big(\det(1-t\phi_{i-1})\det(K_i)\Big)^{(-1)^{i+1}}.
\]
\end{lemma}

\subsection{Geometric interpretation}
\label{sec:interpret}

We will now interpret the factors on the right side of
Lemma~\ref{lem:compute} in terms of Morse theory.

\begin{notation}
Suppose $x,y$ are elements of $\Lambda=\Z[V]((t))$, or matrices with
entries in $\Lambda$, which might depend on the choice of
triangulation $T_1$ and cellular approximation $\rho$. We write
\[
x\approx y
\]
if $x-y=O(t^k)$, where $k$ can be made arbitrarily large
by choosing $T_1$ sufficiently fine and $\rho$ sufficiently close to
the identity.
\end{notation}

\begin{lemma}
\label{lem:interpret}
\begin{itemize}
\item[\rm(a)]
$\prod_{i=0}^{n-1}\det(1-t\phi_i)^{(-1)^i}\approx\zeta$.
\item[\rm(b)]
Under the natural identification $D_*\simeq \Cnov_*$, we have
\[
{\bf K}_i\approx\partial^f_i.
\]
\end{itemize}
\end{lemma}

\begin{proof}
(a)\qua
Let $\hat{f}\co \tilde{X}\to\R$ be a lift of $f$, and let
$\tilde{\Sigma}\eqdef \hat{f}^{-1}(0)$.  The downward gradient flow
of $\hat{f}$ induces partially defined return maps
\[
\varphi^k\co \tilde{\Sigma}\to t^k\tilde{\Sigma}.
\]
The definition \eqref{eqn:zetaIntrinsic} of $\zeta$ is equivalent to
\begin{equation}
\label{eqn:zetaFixedPoints}
\zeta=\exp\left(\sum_{k>0,\;g\in V}\op{Fix}(\varphi^k\circ
t^{-k}g^{-1})\frac{gt^k}{k}\right)\in \Z[V]((t))=\Lambda.
\end{equation}
Here $\op{Fix}(s)$ counts fixed points of the equivariant map $s$
modulo covering transformations, with their Lefschetz signs.

Suppose to begin that $\rho=\op{id}$.  By the machinery used to prove the
Lefschetz fixed point theorem in \cite{brown}, for each $k$ we have
\begin{equation}
\label{eqn:fix}
\sum_{g\in V}\op{Fix}(\varphi^k\circ t^{-k}g^{-1})\cdot g =
\sum_{i=0}^{n-1}(-1)^i\op{Tr}(\phi_i^k) \in\Z[V].
\end{equation}
In this case we have
\begin{equation}
\label{eqn:lefschetz}
\zeta=\prod_{i=0}^{n-1}\det(1-t\phi_i)^{(-1)^i}.
\end{equation}
To see this, it is enough to check that the logarithmic derivatives of
both sides are equal, which follows from equations \eqref{eqn:zetaFixedPoints}
and \eqref{eqn:fix}.

In general, let $H\co \Sigma\times[0,1]\to\Sigma$ be the homotopy
from $\op{id}$ to $\rho$.  In \cite{hutchings-lee} we defined a
natural compactification $\overline{\Gamma} \subset
\Sigma\times\Sigma$ of the graph of $\varphi$.  Using this one can
define a compactified graph $\overline{\Gamma}_t^i$ of
$(H(\cdot,t)\circ\varphi)^i)$ in a similar manner.  Now there exists a
positive integer $N$ such that if the cells in $T_1$ are all contained
in balls of radius $\epsilon$, then the homotopy $H$ can be chosen so
that $\op{dist}(H(t,x),x)<N\epsilon$ for all $t\in[0,1]$ and
$x\in\Sigma$. (Cf the construction of $H$ in \cite{spanier}; by
carefully controlling each intermediate step in the homotopy, the above
claim may be achieved.)
Also, the set of fixed points of $\varphi^i$ lies in the interior of 
$\bar{\Gamma}^i$ under the diagonal map $\Sigma\to \Sigma\times\Sigma$ by definition of $\bar{\Gamma}^i$, and is compact as a
consequence of Assumption~\ref{assumption:admissible}. It follows that
for any positive integer $k$ we can choose $\epsilon$ so that for all
$i\le k$ and all $t\in[0,1]$, the compactified graph
$\overline{\Gamma}_t^i$ does not cross the diagonal in
$\Sigma\times\Sigma$.  Then equation \eqref{eqn:fix} will hold up to
order $k$, and therefore so will equation \eqref{eqn:lefschetz}.

(b)\qua
If $p$ is a critical point of index $i$, then
\begin{equation}
\label{eqn:bdp=dbp}
\partial[\mc{D}(p)]=[\mc{D}(\partial^fp)]
\end{equation}
where the brackets indicate the fundamental class of the
compactification of the descending manifold \cite{hutchings-lee}.  Now
suppose again that $\rho=\op{id}$. Recall from equation
\eqref{eqn:bfK} that the matrix ${\bf K}_i$ sends
$\widehat{\mc{D}_0(p)}$ to a linear combination of cells of the form
$\widehat{\mc{D}_0(q)}$, where $q$ is a critical point of index $i-1$.
(We will henceforth omit the hats $\hat{}$ on $\mc{D}_0$ or $\mc{F}$
when $\rho=\op{id}$, since in this case hatted and unhatted versions
can be identified.)  In fact,
\[
{\bf K}_i ({\cal D}_0(p))={\cal D}_0 (\partial^f (p)).
\] 
To see this, note that from the
definition of ${\bf M}_i$ in equation
\eqref{eqn:bigMatrix}, we have
\[
[\mc{D}(p)] = \mc{D}_0(p)+\sum_{k=0}^\infty
t^{k+1}\mc{F}(\phi^k{\bf M}_i(\mc{D}_0(p))).
\]
(Here the initial descending manifolds $\mc{D}_0$ and initial downward
flow $\mc{F}$ are defined as in \S\ref{sec:complex}, but using
$\tilde{X}$ and $\tilde{\Sigma}$ instead of $X$ and $\Sigma$.)
Applying equation \eqref{eqn:bigMatrix} to this gives
\begin{equation}
\label{eqn:bdp}
\partial[\mc{D}(p)]= {\bf K}_i(\mc{D}_0(p))+\mbox{(terms without
initial descending manifolds).}
\end{equation}
Equations \eqref{eqn:bdp=dbp} and \eqref{eqn:bdp} imply that
$\partial_i^f={\bf K}_i$ when $D_*$ (the domain/range of ${\bf K}_*$) is
identified with $CN_*$.

The case $\rho\neq\op{id}$ can be handled similarly to part (a).
\end{proof}

We can now complete the proof of Theorem~\ref{thm:main}.
Lemmas~\ref{lem:interpret}(b) and \ref{lem:computeTorsion} imply that
\[
\prod_{i=1}^n\det(K_i))^{(-1)^{i+1}}\approx\tau(\Cnov_*)(\xi).
\]
Together with Lemmas~\ref{lem:technical},~\ref{lem:compute}, and
\ref{lem:interpret}(a), this implies that Theorem~\ref{thm:main} holds
up to order $k$ for all $k$.  \hfill \qedsymbol

\section{The 3--dimensional case and Seiberg--Witten theory}
\label{sec:sw}

We now review from \cite{hutchings-lee} the definition of the
Morse-theoretic invariant $I_3$, and the background for
Conjecture~\ref{conj:main} relating this invariant to Seiberg--Witten
theory.  We will then prove Theorem~\ref{thm:sw}, relating this
invariant to Turaev torsion.

\subsection{Motivation from Seiberg--Witten theory}

Let $X$ be a closed connected oriented smooth 3--manifold with
$b_1(X)>0$.  Let $\mathfrak{s}$ be a spin-c structure on $X$.  This
determines a $U(2)$--bundle $S\to X$ with a Clifford action of $TX$ on
$S$.  A section $\psi$ of $S$ and a connection $A$ on $\det(S)$
satisfy the Seiberg--Witten equations with perturbation $\omega$ if, in
the notation of \cite{kronheimer-mrowka},
\[
\begin{split}
\not\partial_A\psi&=0,\\
\rho(F_A)&=i\sigma(\psi,\psi)+i\rho(\omega).
\end{split}
\]
The Seiberg--Witten invariant $\SW(\mathfrak{s})$ counts solutions to these
equations modulo gauge equivalence.  (For more on 3--dimensional
Seiberg--Witten invariants see
eg \cite{lim,meng-taubes,okonek-teleman}.)

Let us choose the perturbation to be $\omega=r* df$, where $f\co X\to
S^1$ is harmonic, $*$ denotes the Hodge star, and $r$ is a real
number.  By perturbing the metric, we may arrange that $f$ is a Morse
function.  Away from the critical points, the spinor bundle $S$ splits
into eigenspaces of Clifford multiplication by $df$,
\begin{equation}
\label{eqn:clifford}
S=E\oplus(E\otimes K^{-1}),
\end{equation}
where $K^{-1}\eqdef\Ker(df\co TX\to\R)$.

Taking $r\to\infty$, one expects that for a Seiberg--Witten solution
the zero set of the $E$ component of
$\psi$ to become parallel to $\nabla f$.  (The energy of the
Seiberg--Witten solution will be concentrated along this zero set.  For
detailed analysis see \cite{taubes:sw=gr} and its sequels.)  This suggests that
$\SW(\mathfrak{s})$ counts unions of closed orbits and flow lines of
$\nabla f$ starting and ending at critical points, 
whose total homology class is Poincar\'{e} dual to
$c_1(E)$.

The above homological condition implies that in our union of closed
orbits and flow lines there is precisely one flow line starting at
each index 2 critical point and ending at each index 1 critical point.
(See \cite{hutchings-lee}.  There are no index 0 or 3 critical points
because $f$ is harmonic, and there are equally many index 1 and index
2 points because $\chi(X)=0$.)  In other words, in the notation of
\S\ref{sec:euler}, our union of closed orbits and flow lines lives in
$H_1(X,v)$, where $v=-\nabla f$.  As in Remark~\ref{remark:torus}, the
counting of closed orbits is related to Taubes' counting of
pseudoholomorphic tori in symplectic 4--manifolds
\cite{taubes:counting}, which indicates that we should allow closed
orbits to be multiply covered when they are elliptic, but not when
they are hyperbolic.

\subsection{The definition of $I_3$}

We now want to define $I_3(\mathfrak{s})$ to be a signed count of such
unions of closed orbits and flow lines.  A convenient way to do so is
to use generating functions as follows.  Choose orderings of the index
1 and index 2 critical points.  Let $\mc{P}^{ij}$ denote the set of
flow lines from the $i^{th}$ index 2 point to the $j^{th}$ index 1
point.  Define the path matrix $P$ by
\[
P^{ij}\eqdef\sum_{\gamma\in\mc{P}^{ij}}\epsilon(\gamma)[\gamma].
\]
Here $[\gamma]\in H_1(X,v^{-1}(0))$ is the homology class of $\gamma$
(oriented downward), and $\epsilon(\gamma)$ is the sign of $\gamma$ as
in \S\ref{sec:morse}.  The entries of $P$ live in the Novikov ring of
the relative homology group $H_1(X,v^{-1}(0))$ with grading given by
minus intersection number with $\Sigma$.

Note that $\det(P)$, regarded as a $\Z$--valued function on
$H_1(X,v^{-1}(0))$, is supported on $H_1(X,v)$.  Also, the subset of
$\Nov(H_1(X,v^{-1}(0)))$ consisting of functions supported on
$H_1(X,v)$ is a $\Lambda$--submodule.  A generating function counting
unions of closed orbits and flow lines of the type we want is now
given by
\[
I_3^v\eqdef\zeta\cdot\det(P)\in\Nov(H_1(X,v)).
\]
In the above equation, `$\cdot$' denotes the $\Lambda$ action. The
closed orbits are counted correctly as a result of the product formula
for the zeta function \eqref{eqn:productFormula}.

Let $j_v\co \Spinc(X)\to H_1(X,v)$ denote the map that sends a spin-c
structure to the Poincar\'{e}--Lefschetz dual of $c_1(E)$, where $E$ is
the line bundle defined in \eqref{eqn:clifford}.  The map $j_v$ is an
$H_1$--equivariant isomorphism.

\begin{definition}
\cite{hutchings-lee}\qua Regarding $I_3^v$ as a function $H_1(X,v)\to \Z$,
we define
\[
I_3\eqdef I_3^v\circ j_v\co \Spinc(X)\to\Z.
\]
\end{definition}

This definition makes sense for any Morse function $f\co X\to S^1$ with
no index 0 or 3 critical points, even if $f$ is not harmonic.  The
calculation below will show that $I_3$ does not depend on $f$, except
that there is a global sign ambiguity in $I_3$ due to the orientation
choices we made.  (Also $I_3$ depends on the sign of $\theta\in
H^1(X;\Z)$ when $b_1(X)=1$.)

\subsection{Relation with Turaev torsion}
\label{sec:conclusion}

To relate $I_3$ to torsion, we note that the isomorphism $i_v\circ
j_v\co \Spinc(X)\to\Eul(X)$ does not depend on $v$.  It follows that there is a
canonical isomorphism $\iota\co \Spinc(X)\to\Eul(X)$.  (This isomorphism
was first defined by Turaev \cite{turaev:euler} in a different but equivalent
way.  The inverse map sends a smooth Euler structure represented by a
nonsingular vector field $u$ to the spin-c structure whose spin bundle
is $\C u\oplus u^\perp$ with a standard Clifford action.)  In summary,
we have the following commutative triangle:
$$\begin{array}{ccc}
H_1(X,v)\;\;\;\;\; & \stackrel{i_v}{\longrightarrow} &
\;\;\;\;\;\mbox{Eul}(X) \\
\;\;\;\;\;\mbox{\small $j_v$} \mbox{\large $\nwarrow$} & & \mbox{\large
$\nearrow$} \mbox{\small $\iota$}\;\;\;\;\; \\
& \Spinc(X) &
\end{array}$$

\begin{proof}[Proof of Theorem~\ref{thm:sw}]
It is enough to show that
\[
I_3^v\circ i_v^{-1}=\pm T(X;o)\co \Eul(X)\to\Z.
\]
By Theorem~\ref{thm:main} and the definition of $T$, this is equivalent
to asserting that for some orientation choice, we have
\[
\det(P)(i_v^{-1}(\xi))=\tau(\Cnov_*)(\xi)(0)
\]
for all Euler structures $\xi$.  If $\xi_0$ is a reference Euler
structure and $\gamma\eqdef i_v^{-1}(\xi_0)\in H_1(X,v)$, then this
equation is equivalent to
\[
\det(P)(\cdot+\gamma)=\tau(\Cnov_*)(\xi_0)(\cdot)\co H_1\to\Z.
\]
This last equation follows from the definition of $\tau(\Cnov_*)$ and
Example~\ref{example:determinant}.
\end{proof}

%
%

\end{document}